\documentclass[graybox]{svmult}

\usepackage{mathptmx}       
\usepackage{helvet}         
\usepackage{courier}        
\usepackage{type1cm}        
\usepackage{makeidx}         
\usepackage{graphicx}        
\usepackage{multicol}        
\usepackage[bottom]{footmisc}
\usepackage[T1]{fontenc}
\usepackage{amsmath}
\usepackage{amsfonts}
\usepackage[polutonikogreek, english]{babel}

\newcommand{\R}{{\mathbb{R}}}

\newcommand{\N}{{\mathbb{N}}}

\newcommand{\F}{{\mathbb{F}}}
\newcommand{\A}{{\mathbb{A}}}

\newcommand{\pv}{\par\vspace{1ex}}

\usepackage{graphicx}

\makeindex             

\begin{document}

\large

\title*{Euler's Series for Sine and Cosine. An Interpretation in Nonstandard Analysis}
\author{Piotr B{\l}aszczyk and Anna Petiurenko}
\institute{Piotr B{\l}aszczyk \at Pedagogical University of Krak\'{o}w, ul. Podchor\k{a}\.{z}ych 2, 30-084 Krak\'{o}w, Polska \email{\mbox{piotr.blaszczyk@up.krakow.pl}}
\and Anna Petiurenko \at Pedagogical University of Krak\'{o}w, ul. Podchor\k{a}\.{z}ych 2, 30-084 Krak\'{o}w, Polska \email{anna.petiurenko@up.krakow.pl}}
%
%
\maketitle

\abstract*{In chapter VIII of \textit{Introductio in analysin infinitorum}, Euler derives a series for sine, cosine, and the formula $e^{iv}=\cos v+i\sin v$  His arguments employ infinitesimal and infinitely large numbers and some strange equalities. We interpret these seemingly inconsistent objects within the field of hyperreal numbers. We show that any non-Archimedean field provides a framework for such an interpretation. Yet, there is one implicit lemma underlying Euler's proof, which requires specific techniques of non-standard analysis. Analyzing chapter III of \textit{Institutiones calculi differentialis} reveals Euler's appeal to the rules of an ordered field which includes infinitesimals -- the same ones he applies deriving series for $\sin v$, $\cos v$, and $e^{v}$. 
}

\abstract{In chapter VIII of \textit{Introductio in analysin infinitorum}, Euler derives a series for sine, cosine, and the formula $e^{iv}=\cos v+i\sin v$. His arguments employ infinitesimal and infinitely large numbers and some strange equalities. We interpret these seemingly inconsistent objects within the field of hyperreal numbers. We show that any non-Archimedean field provides a framework for such an interpretation. Yet, there is one implicit lemma underlying Euler's proof, which requires specific techniques of non-standard analysis. Analyzing chapter III of \textit{Institutiones calculi differentialis} reveals Euler's appeal to the rules of an ordered field which includes infinitesimals -- the same ones he applies deriving series for $\sin v$, $\cos v$, and $e^{v}$. 
}

\section{Introduction}
\label{sec:1}
In chapter VII of \textit{Introductio in analysin infinitorum} (\S\S 114--125), Euler introduces the famous number $e$ and derives a series for $e^z$. Let us focus on the technique applied, especially infinitesimals and infinitely large numbers that he refers to in his arguments. 

The first paragraph of chapter VII reads:

\begin{quote} Since $a^0 = 1$, when the exponent on $a$ increases, the power itself 
increases, provided $a$ is greater than 1. It follows that if the exponent is 
infinitely small and positive, then the power also exceeds 1 by an infinitely small 
number. Let $\omega$ be an infinitely small number, or a fraction so small that, 
although not equal to zero, still $a^\omega = 1 + \psi$, where $\psi$ is also an infinitely small 
number. From the preceding chapter we know that unless $\psi$ were infinitely small then 
 neither would $\omega$ be infinitely small. It follows that $\psi= \omega$, or $\psi > \omega$, or $\psi<\omega$.
 Which of these is true depends on the value of $a$, which is not now known, so let $\psi=k\omega$. Then we have 
 $a^\omega=1+k\omega$,  and with $a$ as the base 
for the logarithms, we have $\omega=\log (1+k\omega)$. (Euler 1988, 92) \end{quote}

\pv Infinitesimals referred to in this section are involved in the mathematical process on a par with standard quantities. They are subject to the trichotomy law,  $\psi= \omega$, or $\psi > \omega$, or $\psi<\omega$.
They are also terms of operations such as exponents or logarithms. 

Let us proceed further.

\begin{quote} Since  $a^\omega=1+k\omega$, we have $a^{j\omega}=(1+k\omega)^j$,
whatever value we assign to $j$. It follows that 
\[a^{j\omega}= 1+\dfrac 1jk\omega + \dfrac {j(j-1)}{1\cdot 2}k^2\omega^2+\dfrac {j(j-1)(j-2)}{1\cdot 2\cdot 3}k^3\omega^3+\&c.   \]
If now we let $j=\frac{z}{\omega}$, where $z$ denotes any finite number, since $\omega$ is infinitely small, 
then $j$ is infinitely large.
Then we have $\omega=\frac zj$,  where $\omega$ is represented by a  fraction with an infinite denominator, 
so that $\omega$ is infinitely small, as it should. 
When we substitute $\frac zj$
 for  $\omega$ then 

\begin{equation*}\label{E6}
\begin{aligned}
a^{z}=(1+\tfrac{kz}j)^j= 1+\dfrac 11kz + \dfrac {1(j-1)}{1\cdot 2j}k^2z^2+\dfrac {1(j-1)(j-2)}{1\cdot 2j\cdot 3j}k^3z^3+\\
 +\dfrac {1(j-1)(j-2)(j-3)}{1\cdot 2j\cdot 3j\cdot 4j}k^4z^4 + \&c. 
\end{aligned}
\end{equation*}

This equation is true provided an infinitely large number is substituted for $j$, but then $k$ is a 
finite number depending on $a$, as we have just seen. (Euler 1988, 93).\footnote{
In Blanton's translation (Euler 1988), instead of the original sign $\&c.$, there are three dots $\ldots$\,. In section 4.3 below, we comment on that.}\end{quote}

\pv Next to infinitesimals and infinite numbers, the above section introduces finite numbers. It also provides some relationships: the product $\frac z\omega$ is infinite, given $z$ is finite and $\omega$ infinitesimal, or 
$\frac zj$ is infinitesimal, given $j$ is infinite. Moreover, it considers sums with infinitely many terms.
Note also that Euler processes fractions involving infinitesimals and infinite numbers as if were usual fractions: $\omega=\frac zj$, given $j=\frac z\omega$. In the next sections, Euler puts $z=1$. Consequently, the reciprocal of an infinite number is infinitesimal, and of an infinitesimal is an infinite number.

In section 116, we find somewhat strange arithmetic terms involving infinite numbers.

\begin{quote}Since $j$ is infinitely large, $\frac{j-1}{j}=1$, and the larger the number we 
substitute for $j$, the closer the value of the fraction 
 comes to $1$.  Therefore, if $j$ is a number larger than any assignable number, then 
 $\frac{j-1}{j}$ is equal to $1$. For the same reason  $\frac{j-2}{j}=1$, $\frac{j-3}{j}=1$,  and so forth.
 It follows that $\frac{j-2}{2j}=\frac 12$, $\frac{j-3}{3j}=\frac 13$, $\frac{j-4}{4j}=\frac 14$, 
 and so forth. (Euler 1988, 93--94)\end{quote}
 
 \pv In this paragraph, infinite numbers, i.e., those greater than \textit{any assignable number}, satisfy  a mystifying equality
  $\frac{j-1}{j}=1$. Moreover, on the one hand $\frac{j-1}{j}$ equals $1$, on the other, 
  $\frac{j_1-1}{j_1}$
  gets closer to $1$ than  $\frac{j_2-1}{j_2}$, given $j_1>j_2$. Clearly, the equality sign does not stand here for the strict equality.  Since $\tfrac 1j$ is infinitesimal, $\frac{j-1}{j}$ or $1-\tfrac 1j$ is infinitely close to $1$ rather than equal to it.

  By substituting $\frac 1k$ for $\frac{j-k}{kj}$ in the series for $a^z$, Euler gets
  \[a^{z}=1+ \dfrac{kz}1 + \dfrac {k^2z^2}{1\cdot 2}+\dfrac {k^3z^3}{1\cdot 2\cdot 3}+ 
  \dfrac {k^4z^4}{1\cdot 2\cdot 3\cdot 4} + \&c.   \]

Assuming $z=1$, he finds a series for $a$, 
 \[a=1+ \dfrac{k}1 + \dfrac {k^2}{1\cdot 2}+\dfrac {k^3}{1\cdot 2\cdot 3}+ 
  \dfrac {k^4}{1\cdot 2\cdot 3\cdot 4} + \&c.   \]

 Taking $k=1$, he can
define number $e$ as follows
\[e= 1+\frac 11+\frac 1{1\cdot 2}+\frac 1{1\cdot 2\cdot 3}+\frac 1{1\cdot 2\cdot 3\cdot 4}+\&c.\  \mbox{in\ infinitum}.\] 

Going through that argument once again, this time with the number $e$ instead of $a$, we get the following:
\begin{align}\label{E0}
&e^\omega=1+\omega, \nonumber\\
&e^{j\omega}=(1+\omega)^j, \nonumber\\
&e^z=(1+\tfrac zj)^j,\ \ \ \mbox{where}\ \ j=\dfrac z\omega.
\end{align}

Finally, assuming $\binom j m\tfrac 1{j^m}=\tfrac 1{m!}$, from (\ref{E0}), the series for $e^z$ follows
\[e^{z}= 1+\frac z1+\frac {z^2}{1\cdot 2}+\frac {z^3}{1\cdot 2\cdot 3}+\frac {z^4}{1\cdot 2\cdot 3\cdot 4}+
\&c.\  \mbox{in\ infinitum}.\]

There are intermezzos in Euler's argument for the series expansion of $e^z$ providing numerical examples of some terms, such as this one. 
When $a=10$ and $k\omega=\dfrac 1{1000000}$, then from the tables of logarithm Euler  derives the value of $\omega$ as follows 
$$\log{\Big(1+\dfrac{1}{1000000}}\Big)= \log{\dfrac{1000001}{1000000}} =0.00000043429=\omega.$$ 

As a result, $k=2.30258$. 

These calculations aim  to illustrate that $k$ depends on $a$.
We assume that in this case, $\omega$ is ``a fraction so small" rather than infinitesimal. 

Euler also estimates  the number $e$, namely, it is equal to 
$$2.71828182845904523536028,\ \ \mbox{or}\ \ \ 2.718281828459\&c.$$ 

 In sections VII and then VIII, finite numbers are common numbers usually given in the decimal form. However,  products of infinitesimals and infinite numbers, like  
$j\omega$,   also represent finite numbers. Since Euler does not exemplify infinitesimals and infinitely large, they seem to constitute a technique of proving some results concerning finite numbers. Indeed, the implicit procedure is this: when a finite number $z$ is given, take any infinite number $j$ -- more precisely, a non-standard integer $j$ -- then, $z$ and $j$ determine infinitesimal, $\omega=\omega(z,j)=\dfrac{z}{j}$, and the formula (\ref{E0}) follows. Although the series for $e^z$ consists of infinitely many terms, to estimate $e^z$, Euler considers finitely many of them and does not examine the rest.
 
In chapter VIII, Euler applies the same technique to derive series for sine and cosine starting from the formulas
\[\cos nz= \frac{(\cos z + i \sin z)^n+ (\cos z -i\sin z)^n}{2},\]
 \[\sin nz= \frac{(\cos z + i \sin z)^n- (\cos z -i\sin z)^n}{2i}. \]

He reached these results by reviewing then-current trigonometry through complex numbers.
At the end of that chapter, he combines these results and the formula (\ref{E0}). Once again, a reference to infinitesimals, infinite numbers, and specific understanding of equality substantially contributed to his argument. Section 138 reads:

\begin{quote}let $z$ be an 
infinitely small arc and let $n$ be an infinitely large number $j$, so that $jz$ has a 
finite value $v$. Now we have $nz=v$ and $z=\dfrac vj$, so that $\sin z=\dfrac vj$, and $\cos z=1$. With these substitutions 
\[\cos v=\dfrac{\Big(1+\dfrac {iv} j\Big )^j+\Big(1-\dfrac{iv}j\Big)^j}{2},\ \ \ \mbox{and}\ \ \ 
\sin v =\dfrac{\Big(1+\dfrac {iv} j\Big )^j-\Big(1-\dfrac{iv}j\Big)^j}{2i}.\]

In the preceding chapter we saw   that $(1+z/j)^j=e^z$ 
where $e$ is the base of the natural logarithms. When we let $z=iv$ and then $z=-iv$ we obtain
\[\cos v=\dfrac{e^{iv}+e^{-iv}}2,\ \ \ \mbox{and}\ \ \ \sin v=\dfrac{e^{iv}-e^{-iv}}{2i}.\]
From these equations we understand how complex exponentials can be expressed by real sines and cosines, since 
$e^{iv}=\cos v+i\sin v$, and $e^{-iv}=\cos v -i\sin v$. (Euler 1988, 111--112)\end{quote}

To sum up, Euler's mathematical toolbox includes: (1) the multiplicative inverse of infinitesimal is infinitely large, (2) the multiplicative inverse of  infinitely large is infinitesimal, (3) for a finite number $v$ and  infinite number $j$, there is an  infinitesimal $\omega$   such that $j\omega=v$, (4) infinite numbers $j, m$ satisfy the equality $\tfrac{j-m}{mj}=\tfrac 1m$, and  $\binom jm\tfrac 1{j^m}=\tfrac 1{m!}$, (5) binomial theorem applies to an infinite number $j$, i.e.,
{$(1+x)^j= 1+\binom j 1x+\binom j 2x^2+...$\,.  
 Chapter VIII of \textit{Introductio in analysin infinitorum} brings in another assumptions, namely (6) $\sin z =z$ and $\cos z=1$, for \mbox{infinitesimal $z$,}  (7) a sum of infinitely many  infinitesimals is infinitesimal.

In this paper, we interpret these assumptions within non-standard analysis and focus on Euler's series expansion for sine and cosine developed in chapter VIII  of \textit{Introductio in analysin infinitorum}. In the middle of the 18th century, the number $e$ and series for $e^z$ were brand new topics.  Series for sine and cosine enable to present Euler's novel technique against the ancient old mathematical problem. Accordingly, in section 2, we briefly describe the way Ptolemy and Newton dealt with the sine.    In section 3, we review the basics of Euler trigonometry. Section 4 includes a detailed analysis of how Euler expands sine and cosine into series and what kind of series these are. In section 5, we introduce basic concepts of non-standard analysis and interpret assumptions  (1) through (6) listed above. Thesis (7) requires more advanced techniques of non-standard analysis, and we will not discuss it in this paper.  In section 6,  we support our interpretation by analyzing chapter III of  \textit{Institutiones calculi differentialis}.

\section{Forerunners}
\label{sec:2}
The mathematical concept of sine differed through the ages. In this section, we focus on two vital contributions: Ptolemy's \textit{Table of Chords} and Newton's series of sine.  

\subsection{Ptolemy} Ptolemy identified sine with a chord of a circle considered in respect to the arc it subtends and the respective central angle; see the chord $AB$, arc $AB$ and the angle $AOB$ in Fig.\ref{figP}.   In the \textit{Almagest}, he  tallies  arcs and  chords due to the formula 
$$arc:360::chord:120.$$ 

When an arc is such and such part of the circle (divided into 360 parts),  the chord is such and such part of the diameter (divided into 120 parts).  Ptolemy considered the division of the circle into 720 parts, i.e.,  each $1/2^{\circ}$, and managed to determine the chord for every arc from the $1/2^{\circ}$  to $180^{\circ}$. During that process, he combined Euclid's theory of similar figures and Babylonian arithmetic.\footnote{See: (Ptolemy 1984), Book I, \S\S\, 10--11.} 

Euclid's proportion theory, the foundation of similar figures, presupposed that ratios concern magnitudes of the same kind. Consequently,  in Greek mathematics, ratios such as $arc:chord$ were not legitimate objects. Newton and Euler employed novel techniques, which enabled them to relate arcs and line segments. A profound change in the theory of proportion paved the way to their series for sine. In modern terms, it was that operations in ordered fields replaced the ancient technique of transforming ratios.\footnote{See section 6.1 below.} 
\begin{figure} 
\centering
{\includegraphics[scale=1]{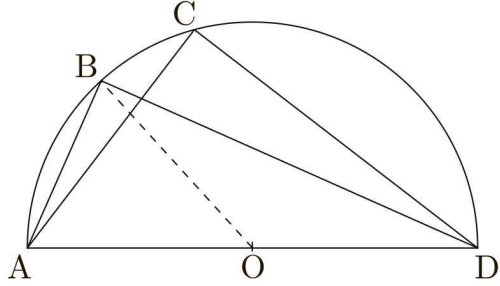}}
\caption{Ptolemy, \textit{Almagest}, p. 51  (the center $O$ and the dashed line added).}
\label{figP}
\end{figure}


 \subsection{Newton} \textit{Analysis by Equations of an Infinite Number of Terms} is another treatise in which Newton demonstrates his  technique of infinite series. We review the results concerning sine and its arc.\footnote{See (Newton 1745), pp. 336--338, (Guicciardini 2009), ch. III.} In Fig.\ref{figN1}, line $TB$ is the tangent to the semicircle  $ADLE$ and $AD$ is the sine of the angle $ACD$. Newton seeks to determine the arc $AD$ in terms of $AB$, given $AB=x$. To this end, he introduces the ``indefinitely small rectangle" $BGHK$, and sets $AC=1/2$. The point $H$   lies on $TD$ at the same time it is a vertex of the triangle $HDG$ and lies on the semicircle. Following Descartes' technique,\footnote{On how Descartes interpreted the Pythagorean theorem see (B{\l}aszczyk, Mr{\'o}wka, Petiurenko 2020, \S\S\,6--7).} Newton finds that $BD$ equals $\sqrt{x-x^2}$. From the similarity of triangles $\triangle TDB$, 
$\triangle HDG$ and $\triangle CDB$, he derives  proportions   
   \[BK:HD::BT:DT,\ \ BT:DT::BD:DC,   \]
and gets  $BK:HD::\sqrt{x-x^2}:1/2$.  Ratio  $HD:KB$ represents  ``the moment of the arc AD" to ``the moment of the base AB".   Once again, due to Descartes' 
technique,\footnote{See (B{\l}aszczyk 2022).} the proportion $HD:BK::DC:BD$ is turned into the subsequent equality
 \[\frac{HD}{BK}=\frac{1}{2\sqrt{x-x^2}}.\]
  
  By the binomial theorem, Newton expands $\frac{HD}{BK}$ into a series,
\[\frac{1}{2\sqrt{x-x^2}}= \frac 12x^{-1/2}+\frac 14 x^{1/2}+\frac 3{16}x^{3/2}+\frac 5{32}x^{5/2}+ \frac{63}{512}x^{9/2},\,\&c.\]

Via term-wise integration, he gets that ``the length of the arc is"

\[x^{1/2}(1+\frac 16 x+\frac 3{40}x^2+\frac 5{112}x^{3}+ \frac{35}{1152}x^4+\frac{63}{2816}x^5,\,\&c.).\]

Then,  Newton  continues: ``After the same manner by supposing CB to be $x$, the radius CA to be 1, you will find the arc LD to be"
\begin{equation}\label{Narc}x+\frac 16x^3  +\frac 3{40}x^5+\frac 5{112}x^7+\&c.\end{equation} 
 
 That is, given $BC=\sin \angle DCL=x$, $CA=1$, the arc $DL$ equals  the series (\ref{Narc}). In  modern terms
 \[\arcsin x=x+\frac 16x^3  +\frac 3{40}x^5+\frac 5{112}x^7+...\,.\]

\begin{figure} 
\centering
{\includegraphics[scale=1]{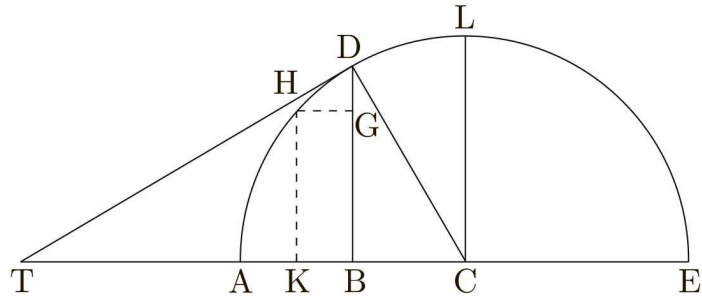}}
\caption{Newton, \textit{Analysis of Equations}, p. 336.}
\label{figN1}
\end{figure}

In the next section, Newton seeks to determine the sine in terms of its arc. In Fig.\ref{figN2}, the line  
$AB$ is the sine of the angle $\alpha A D$.  Given $\alpha D=z$, $AB=x$, $A\alpha=1$, due to (\ref{Narc}), the equality obtains
\[z=x+\frac 1 6 x^3+\frac{3}{40}x^5+\frac{5}{112}x^7,\&c.\]

By his ingenious technique of finding inverse series,  Newton determines the relation between $z$ and $x$ as follows
\[x=z-\frac 16 z^3+ \frac 1{120}z^5-\frac 1{5040}z^7,\,\&c.\]

In modern terms, it is the series for sine, namely 
\[AB=sin z=z-\frac 1{3!}z^3+\frac 1{5!}z^5-\frac 1{7!}z^7+...\,.\]

\begin{figure} 
\centering
\includegraphics[scale=1.2]{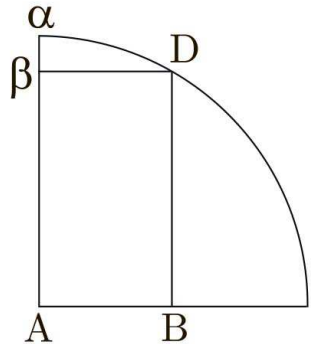}
\caption{Newton, \textit{Analysis of Equations}, p. 338.}
\label{figN2}
\end{figure}

\subsection{Ptolemy--Newton--Euler} Newton's analysis of sine involves tangent to a circle. As a result, he considers the sine of half of the respective angle -- of the say Ptolemy angle. Euler, as we will see, settles  Newton's sine and the tangent line on the unit circle in a way a modern reader takes as normal. At the same time, the tangent gets a new meaning: by definition, it is a ratio $\frac {\sin z}{\cos z}$. 

Since both Newton and Euler employ infinitesimals,  their techniques are commonly considered intuitive versions of modern calculus. While Newton's infinitesimals rely  on intuition, 
Euler defines these seemingly strange objects. In chapter VIII of (Euler 1748), he applies them to derive the sine and cosine series. In the following sections,  we show how to interpret Euler proof in a  branch of calculus, which instead of the concept of limit, applies infinitesimals, that is, in non-standard analysis.\footnote{Regarding the idea of calculus without the concept of limit see (Błaszczyk, Major 2014).}

\pv
\section{Setting the stage. Trigonometry}
\label{sec:3}

Euler's \textit{Introductio in analysin infinitorum} (Euler 1748) is considered the beginning of modern trigonometry.
Here is the founding supposition of its chapter VIII (\S\S\,126--142):
\begin{quote}We let the radius, or total sine, of a circle be equal to 1, then it is clear 
enough that the circumference of the circle cannot be expressed exactly as a 
rational number. An approximation of half of the circumference of this circle is 

 \pv 3.141592653589793238462643383279502884197169399375105820974944592
 
 3078164062862089986280348253421170679821480865132723066470938446+. 

\pv For the sake of brevity we will use the symbol $\pi$ for this number. We say, then, 
that half of the circumference of a unit circle is  $\pi$, or that the length of an arc of 
180 degrees is  $\pi$. (Euler 1988, 101)\end{quote}

The phrase \textit{total sine} refers to the diameter of the circle. Euler's account is novel in many respects. First of all, the unit circle sets up a new standard.

\begin{quote}We always assume that the radius of the circle is 1 and let $z$ be an arc 
of this circle. We are especially interested in the sine and cosine of this arc $z$. 
Henceforth we will signify the sine of the arc $z$ by $\sin z$. Likewise, for the cosine 
of the arc z we will write $\cos z$.  Since $\pi$ is an arc of 180$^{\circ}$, $\sin 0\pi=0$ and $\cos 0\pi=1$. Also
$\sin \pi/2=1$, $\cos \pi/2=0$, $\sin \pi=0$, $\cos \pi=-1$, $\sin 3\pi/2=-1$, $\cos 3\pi/2=0$, 
$\sin 2\pi=0$, and $\cos 2\pi=1$. Every sine and 
cosine lies between +1 and -1. Further, 
we have $\cos z=\sin(\pi/2-z)$,  $\sin z=\cos (\pi/2-z)$. We also have $\sin^2 z+\cos^2 z=1$.
Besides these 
notations we mention also that $\tan z$ indicates the tangent of the arc $z$, $\cot z$ 
for the cotangent of arc $z$. We agree that $\tan z=\frac{\sin z}{\cos z}$ and
$\cot z=\frac{\cos z}{\sin z}=\frac 1 {\tan z}$,
all of which is known from trigonometry. (Euler 1988, 102)\end{quote}

\begin{figure}[h!]
\begin{center}
\includegraphics[scale=0.9]{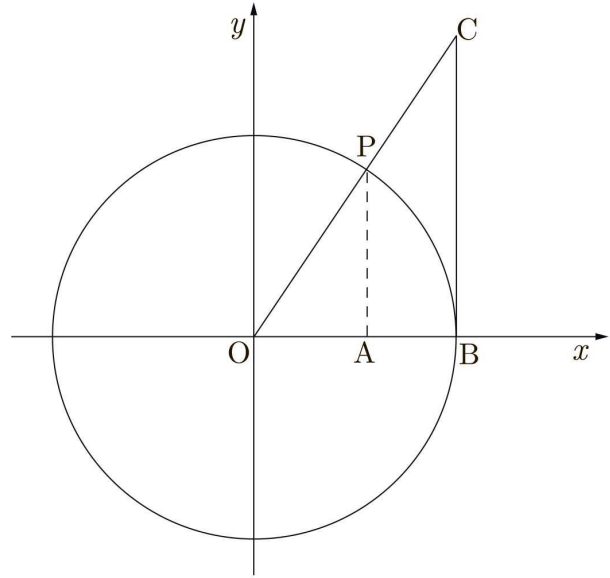}
\caption{Unit circle} \label{figSin}
\end{center}
\end{figure}


\pv Through subsequent sections, Euler surveys standard identities like the sine and cosine sum and difference laws, or half-angle formulas. In Fig.\ref{figSin}, we represent his model of the unit circle with sine, cosine, and tangent. Although (Euler 1748) does not include such a diagram nor does Euler define sine and cosine, it is implicitly hidden behind his considerations. Note that in Fig.\ref{figSin}, sine is parallel to tangent. Thus, in Fig.\ref{figN1}, it would be a line parallel to $TD$, or perpendicular to $DC$,  passing through $A$, and equal to $DC$.

In section 133, Euler introduces another revolutionary trick.   Starting with the Pythagorean theorem in trigonometric stylization, 
\begin{equation}\label{E1}\sin^2 z+\cos^2 z=1,\end{equation}
he factors it  in the field of complex numbers, and gets the following equality
\begin{equation}\label{E2}(\cos z+ i \sin z)(\cos z-i \sin z )=1.\end{equation}

This seemingly simple motion aims to justify the appeal to complex numbers  in trigonometric studies. Euler comments on it as follows: ``Although these factors are complex, still 
they are quite useful in combining and multiplying arcs".\footnote{Rather than the modern sign $i$, (Euler 1748) uses $\sqrt{-1}$.}

By that technique, he derives the formula currently named after de Moivre.  Comparing the real and imaginary parts of products of complex numbers given in trigonometric forms, he shows that
\[(\cos x\pm i \sin x)(\cos y\pm i \sin y)=\cos (x+y)\pm i \sin (x+y),\]
and
{\normalsize \[(\cos x\pm i \sin x)(\cos y\pm i \sin y)(\cos z\pm i \sin z)=\cos (x+y+z)\pm i \sin (x+y+z).\]}

On these grounds, he reaches  the general conclusion
\begin{equation}\label{E3}(\cos z\pm i\sin z)^n= (\cos nz\pm i\sin nz).\end{equation}   

\pv Glen Van Brummelen, an expert in the history of trigonometry, comments on the introductory  sections of the chapter VIII as follows: ``With the adoption of the unit circle, for the first time the sine and cosine are considered to be ratios of line segments rather than their lengths.
Indeed, at the end of the following paragraph, $\tan z$ and $\cot z$ are
introduced directly as $\frac{\sin z}{\cos z}$ and $\frac{\cos z}{\sin z}$ respectively" (Van Brummelen 2021, 166)

Indeed, for Ptolemy and Newton, sine was a specific line segment. Euler explicitly shows that \textit{his} sine and cosine take negative values. Therefore, they should be considered coordinates of a point on the unit circle rather than ratios of line segments. On other occasions, Euler identified ratios and numbers, thus the actual meaning of the term \textit{ratio}  in the context of his trigonometry is the quotient.

Further, Van Brummelen adds: ``Euler does not define the sine and cosine at all in this text; the primitive
notions remain geometrical". Whether on purpose or not, the absence of definition enables Euler to employ various techniques like Euclidean geometry, analysis in terms of Cartesian coordinates, and complex numbers. 
Under the line segment interpretation of 
$\sin z$ and $\cos z$, the equality (\ref{E1}) means the Pythagorean theorem, where the number 1 in Fig.\ref{figSin} stands for the line    $OP$. One can also interpret complex numbers in (\ref{E2}) geometrically. Then, in Fig.\ref{figSin}, the number 1 stands for the line segment $OB$.   Although the processing of the formula (\ref{E1}) into (\ref{E2}) is correct,  the geometrical content changes during that process.  Euler's trigonometry takes advantage of such equivocations.

\section{The crucial move}

Euler  expands  sine and cosine into series in a seemingly self-evident way  through   a few lines of sections 133--134. 
To show his argument, let us write down (\ref{E3}) in the plus and minus form separately
\begin{equation}\label{E4}(\cos z+ i\sin z)^n= \cos nz + i\sin nz, \ (\cos z - i\sin z)^n= \cos nz - i\sin nz.\end{equation}  

 He explicitly considers formulas that follow from (\ref{E4}), namely\footnote{There is a typo in Blanton's translation (Euler 1988): the formula for $\sin nz$ is missing $i$ in the denominator.}
\begin{equation}
\begin{aligned}
\label{E5}\cos nz= \frac{(\cos z + i \sin z)^n+ (\cos z -i\sin z)^n}{2},\\
 \sin nz= \frac{(\cos z + i \sin z)^n- (\cos z -i\sin z)^n}{2i}. 
\end{aligned}
\end{equation}	

Then he writes: "Expanding the binomials, we obtain the following series":\footnote{1) Euler introduces binomials $\binom nk$  in the explicit fractional form.  We employ modern notation, to get a more compact formula. 2) Instead of  Euler's original sign $\&c.$ ending the series, Blanton's translation (Euler 1988) applies three dots $...$, i.e., the symbol marking infinite series in real analysis.  As we will see in section 4.3 below, Euler's infinite series do have the last term.}

\begin{equation}\label{E6}
\begin{aligned}
\cos nz = \cos^n z -  \binom{n}{2}\cos^{n-2}z \sin^2 z+  \binom{n}{4}\cos^{n-4}z \sin^4 z-\\ -\binom{n}{6}\cos^{n-6}z \sin^6 z +  \&c.
\end{aligned}
\end{equation}
\begin{equation}\label{E7}
\normalsize
\begin{aligned}
\sin nz = \binom{n}{1}\cos^{n-1}z\sin z - \binom{n}{3}\cos^{n-3}z \sin^3z+\binom{n}{5}\cos^{n-5}z\sin^5z-\&c.
\end{aligned}
\end{equation}

Section 134 opens with a passage that encourages our non-standard analysis interpretation.
\begin{quote}Let the arc z be infinitely small, then $\sin z = z$ and $\cos z = 1$. If $n$ is 
an infinitely large number, so that $nz$ is a finite number, say $nz = v$, then, since 
$\sin z=z=\frac vn$, we have
\[\cos v =  1 - \dfrac{v^2}{1 \cdot 2} + \dfrac{v^4}{1 \cdot 2 \cdot 3 \cdot 4} - \dfrac{v^6}{1 \cdot 2 \cdot 3 \cdot 4 \cdot 5 \cdot 6} + \&c,\]
\[\sin v =  v - \dfrac{v^3}{1 \cdot 2\cdot 3} + \dfrac{v^5}{1 \cdot 2 \cdot 3 \cdot 4\cdot 5} - \dfrac{v^7}{1 \cdot 2 \cdot 3 \cdot 4 \cdot 5 \cdot 6\cdot 7} + \&c.\]
It follows that if $v$ 
is a given arc, by means of these series, the sine and cosine can be found. (Euler 1988, 107)\end{quote} 

In the series (\ref{E6}), (\ref{E7}), Euler replaces  every term  $\binom{n}{k}\cos^{n-k}z \sin^kz$
 by $\frac{v^k}{k!}$, given $\cos z = 1$, and $\sin z=\frac vn$. In that process, he assumes the following equalities 
\begin{equation}\label{E8}\binom{n}{k}\cos^{n-k}z \sin^kz =\binom{n}{k}1 \frac{v^k}{n^k} =
    \frac{v^k}{k!}.\end{equation}

In the next sections, we go  through this proof in detail.

\subsection{From finite to infinite} Formulas (\ref{E6}) and (\ref{E7}) follow from (\ref{E4}) and (\ref{E5})
due to the binomial theorem. Thus, at first,  $n$ is supposed to be a finite number. In section 134, Euler explicitly considers $n$ infinitely large.  In non-standard analysis, formulas (\ref{E6}) and (\ref{E7}) are valid whether $n$ is a finite or infinite number, specifically  terms $\binom n k$ make sense for finite and infinite $n, k$. More precisely, when the finite case is valid, the infinite is valid too.
   The same applies to the binomial theorem: one can apply it to terms $(\cos z\pm i\sin z)^n$
   whether $n$ is finite or infinite.

  Allowing infinite $n$ in formulas (\ref{E6}) and (\ref{E7}), we need to interpret sums made of infinitely many terms. Before that, we have to comment on infinitesimals.

\subsection{Infinite and infinitesimal numbers} In (Euler 1748), Euler expressly refers to infinitesimals.  On other occasions, for example, in (Euler 1755), he calls them zeros. Then, he introduces two ways of comparing them: arithmetic and geometric.  The former means a difference and allows interpretation in terms of a relation \textit{is infinitely close},  the latter is a quotient.  Non-standard analysis provides an obvious explanation for the  relation  \textit{is infinitely close}, namely $x\approx y$, meaning $x-y$ is infinitesimal.  Thus, assuming $x, y$ are infinitesimals, the difference $x-y$ is still infinitesimal.  Quotient $\dfrac xy$, however, can be infinitesimal, finite, or infinite.

Indeed, in (Euler 1755),  the inverse of an infinitesimal proves to be an infinite number and vice versa.  Euler also shows that for each infinitesimal there is an infinite number such that their product is finite. He also expressly claims that there are infinitesimals and infinite numbers such that their products are infinite. Assumptions we have already found in chapter VII of (Euler 1748),  Euler (Euler 1755) discusses in detail.

Therefore we interpret the passage opening section 134 above ``Let the arc z be infinitely small, then $\sin z = z$ and $\cos z = 1$. If $n$ is an infinitely large number, so that $nz$ is a finite number, say $nz = v$", as follows,
\[z\approx 0 \Rightarrow \sin z\approx z,\ \cos z\approx 1,\]
 and for an infinitesimal $z$ there is an infinite number $N$  such that $Nz$ is finite.

\subsection{Infinite series vs hyperfinite sums}

Instead of the original term $\&c.$ applied in (Euler 1748), Blanton's translation (Euler 1988) employs three dots sign ending infinite sums,  such as this one
{\normalsize \[\cos nz = \cos^n z -  \binom{n}{2}\cos^{n-2}z \sin^2 z+  \binom{n}{4}\cos^{n-4}z \sin^4 z-\binom{n}{6}\cos^{n-6}z \sin^6 z+...\,.  \]}

It suggests the real analysis interpretation: the limit of the series. However, Euler's infinite sums allow another reading:  these series contain the very last term. To elaborate, we refer to section 107 of  (Euler 1755).
Below we cite Blanton's translation of this section  \textit{in extenso}.

\begin{quote}From this we see that he who would say that when this same series
is continued to infinity, that is,
\[1 + x + x^2 + x^3 +\ldots+ x^{\infty},\]
and that the sum is $1/(1-x)$, then his error would be $x^{\infty+1}/(1-x)$, and
if $x > 1$, then the error is indeed infinite. At the same time, however, this
same argument shows why the series $1+ x+x^2 +x^3 +x^4 +\ldots$, continued
to infinity, has a true sum of $1/(1-x)$, provided that $x$ is a fraction less
than 1. In this case the error $x^{\infty+1}$ is infinitely small and hence equal to
zero, so that it can safely be neglected. Thus if we let $x = \frac 12$, then in truth
\[1+\frac 12+\frac 14+\frac 18+\frac 1{16}+\dots=\frac{1}{1-\frac 12}=2.\]

In a similar way, the rest of the series in which $x$ is a fraction less than 1
will have a true sum in the way we have indicated. (Euler 2000, 60)\end{quote}

\pv Instead of three dots fashion, 
\[1+ x+x^2 +x^3 +x^4 +\ldots\ \ \ \mbox{and}\ \ \ 1+\frac 12+\frac 14+\frac 18+\frac 1{16}+\dots\,,\]
(Euler 1755, \S\,107) presents infinite sums in the following form 
\[1+ x+x^2 +x^3 +x^4 +\&c.\ \ \  \mbox{and}\ \ \ 1+\frac 12+\frac 14+\frac 18+\frac 1{16}+\&c.\]

Thus, the symbol 
\[1+ x+x^2 +x^3 +x^4+\&c.\]
stands for
\[1 + x + x^2 + x^3 +\dots+ x^{\infty},\]
rather than 
\[1 + x + x^2 + x^3 +\ldots\,.\]

 In non-standard analysis, infinite sums containing last terms are legitimate objects and  represent the so called hyperfinite sums. Given $N$ is infinite number and $|a|<1$, the following equation
\[1+a^1+a^2+a^3+\ldots +a^N=\frac{1-a^{N+1}}{1-a},\]
can be easily justified within the non-standard framework. 
Moreover, $\frac{1-a^{N+1}}{1-a}$ is infinitely close to $\frac{1}{1-a}$,
\[\frac{1-a^{N+1}}{1-a}\approx\frac{1}{1-a}.\]

Indeed, the absolute value of the difference $\frac{1-a^{N+1}}{1-a}$ and $\frac{1}{1-a}$ equals $\frac{|a|^{N+1}}{1-a}$.
Since $|a|<1$, the number $|a|^{N+1}$ is infinitesimal. Therefore the product  $\frac{|a|^{N+1}}{1-a}$ is also infinitesimal.\footnote{See rule (\ref{Omega2}) below.}

Euler's original formula 
\[1+\frac 12+\frac 14+\frac 18+\frac 1{16}+\&c.=\frac{1}{1-\frac 12}=2,\]
finds an obvious interpretation in non-standard analysis, namely
\[1+\frac 12+\frac 14+\frac 18+\frac 1{16}+\dots+ \frac 1{2^N}=\frac{1-(\frac 12)^{N+1}}{1-\frac 12}\approx \frac{1}{1-\frac 12}=2.\]

Under the hyperfinite sum interpretation, our only intervention in Euler's text concerns the infinity sign: instead of $\infty$ or $\infty+1$, we employ specific infinite numbers (non-standard natural numbers), like $N$, or $N+1$.\footnote{Through the latest study (Ferrero 2022), Ferrero adopts three dots ending series instead of original $\&c.$ and does not consider an alternative interpretation of Euler's infinite series at all.}

\subsection{Extending trigonometric functions} 
Non-standard analysis allows the study of trigonometric functions like 
$\sin^* x,\ \cos^* x$, that extend real   $\sin x$ and $\cos x$. Moreover, these $^*$functions 
obey all the same identities   
as their real counterparts $\sin x$ and $\cos x$. Thus, in our interpretation, formula (\ref{E5}) takes the following form:
\[\cos^* nz = \dfrac{( \cos^* z + i \sin^* z)^n+( \cos^* z - i \sin^* z)^n}{2},\] 
\[\sin^* nz= \frac{(\cos^* z + i \sin^* z)^n- (\cos^* z -i\sin^* z)^n}{2i},\]
 whether  $n$ is a finite or an infinite number. 

The functional interpretation of Euler's sine and cosine is prone to an allegation of an anachronism.
Nevertheless, the $^*$extension of $\sin x$ and $\cos x$ works no matter what is the range of $x$, or whether they are functions or some other objects. 
As long as $\sin x$ and $\cos x$ are given, so are $\sin^* x$ and $\cos^* x$. 

In what follows, we omit that  $^*$ superscript. 

\subsection{The revised proof} In this section, we sum up  our interpretation sketched in sections 4.1 to 4.4. Let $z$ be infinitesimal, that is $z\approx 0$.
Then $\sin z\approx z$ and $\cos z\approx 1$. Let $N$ be an infinite number such that $Nz$ is finite and set $Nz=v$.
First of all, formulas (\ref{E6}), (\ref{E7}) take the form of hyperfinite sums,\footnote{We skip a discussion of whether $N$ is even or odd because these cases do not affect the general picture.}
\begin{equation}
\normalsize
\label{E10}\cos Nz = \cos^N z -  \binom{N}{2}\cos^{N-2}z \sin^2 z+  \binom{N}{4}\cos^{N-4}z \sin^4 z-\ldots + \sin^N z.\end{equation}
\begin{equation}\label{E11}\sin Nz = \binom{N}{1}\cos^{N-1}z\sin z - \binom{N}{3}\cos^{n-3}z \sin^3z+
\ldots+ \sin^Nz.\end{equation}

Then, for (\ref{E8}), we get
\begin{equation}\label{E12}\binom{N}{k}\cos^{N-k}z \sin^kz \approx \frac{v^k}{k!}.\end{equation}

Formula (\ref{E12}) is valid whether $N$ and $k$ are finite or infinite.

From (\ref{E10}), (\ref{E11}) and (\ref{E12}), relationships (\ref{E13}) and (\ref{E14}) follow,
\begin{equation}\label{E13}\cos v \approx  1 - \dfrac{v^2}{2!} + \dfrac{v^4}{4!} - \dfrac{v^6}{6!} + \dots +\dfrac{v^N}{N!} ,\end{equation}
\begin{equation}\label{E14} \sin v \approx  v - \dfrac{v^3}{3!} + \dfrac{v^5}{5!} - \dfrac{v^7}{7!} + \dots +\dfrac{v^N}{N!}.\end{equation}

 To get these results, in (\ref{E10}), (\ref{E11}), we substitute  $\frac{v^k}{k!}$ for $\binom{N}{k}\cos^{N-k}z \sin^kz$. However, these substitutions involve infinitely many terms. There is no obvious  guarantee that 
 \[\cos^N z -  \binom{N}{2}\cos^{N-2}z \sin^2 z+  \binom{N}{4}\cos^{N-4}z \sin^4 z-\ldots + \sin^N z\]
 is infinitely close to 
 $$1 - \dfrac{v^2}{2!} + \dfrac{v^4}{4!} - \dfrac{v^6}{6!} + \dots +\dfrac{v^N}{N!},$$
 even though for every $k$, finite or infinite,  obtains 
 $$\binom{N}{k}\cos^{N-k}z \sin^kz \approx \frac{v^k}{k!}.$$ 
 
 Indeed, a product of infinite and infinitesimal number can be infinite, finite or infinitesimal.
 MacKinzie and Tuckey (MacKinzie, Tuckey 2001) showed that in this case 
the hyperfinite sums in formulas  (\ref{E10}), (\ref{E11}) are infinitely close to hyperfinite sums in (\ref{E13}), (\ref{E14}) respectively.
  Their proof requires sophisticated techniques. In this study, we take it for granted.

Throughout sections 4.1--4.4, we reiterated our belief that concepts Euler applies in chapter VIII like infinitesimals or infinite numbers, can be given strict meaning in the non-standard analysis. In the next section, 
we show how to define in modern terms infinitesimal, infinite number, the relation \textit{is infinity close}, hyperfinite sum, binomial coefficient $\binom N k$ for infinite $N$ and $k$,  $\sin^* x$, and $\cos^*x$. Any non-Archimedean field provides a framework for most of these definitions. In section 6, through analysis of (Euler 1755), we will show that Euler explicitly derives some rules of a non-Archimedean field, even though he does not employ the concept of an ordered field.\footnote{To be clear, only Hilbert (Hilbert 1899) and (Hilbert 1900) introduced the concept of an ordered field.}  Therefore,  although non-standard analysis assumes the theory of real numbers, we do not need to refer to real numbers. There is one exception, however.   The just mentioned lemma of MacKenzie and Tuckey applies specific techniques of non-standard analysis. In this  topic only,  our analysis may be considered anachronistic. 

\section{Non-standard analysis}
 In the brief introduction to the non-standard analysis that follows, we focus on concepts crucial to our interpretation of chapter VIII of (Euler 1748).\footnote{See (Błaszczyk, Major 2014), (Błaszczyk 2016), (Błaszczyk 2021). Section on hyperfinite sums follows (Goldblatt 1998, ch. III).}

\subsection{The basics of ordered field theory}

 A commutative field $(\mathbb{F},+,\cdot,0,1)$ together with a total order $<$ forms an ordered field when field operations are compatible with the order, that is
\[x<y \Rightarrow x+z<y+z,\ \ x<y,\, 0<z \Rightarrow xz<yz.\]

The field of rational numbers is the smallest ordered field, meaning every ordered field includes fractions $\tfrac mn$, for $m, n\in\N$.

In every ordered field, one can define the absolute value
\[|x|=
\left\{\begin{array}{cc}
\ \,x, & \mbox{if\ $x\geq 0$}, \\
-x, & \mbox{if\ $x<0$},
\end{array} \right.\]

and the following subsets of $\mathbb F$:
\begin{eqnarray*}
 \mathbb L &= & \{x: (\exists n\in \mathbb N)(|x|<n)\},\\
 \mathbb A &=&   \{x: (\exists n\in \mathbb N)(\tfrac 1 n<|x|<n)\},\\
 \Psi             &=& \{x: (\forall n\in \mathbb N)(|x|>n)\},\\
 \Omega      &=& \{x: (\forall n\in \mathbb N)(|x|<\tfrac 1n)\}.
 \end{eqnarray*}

The elements of these sets we call finite,  assignable,   infinitely large, and infinitely small (infinitesimals)  numbers, respectively.
Infinitesimals and infinite numbers can be also defined via   assignable numbers,
\begin{equation} \Omega=\{x: (\forall a\in\A)(|x|<|a|)\}, \end{equation}
\begin{equation} \Psi= \{x: (\forall a\in\A)(|x|>|a|)\}. \end{equation}

Here are some obvious relationships between these kinds of elements, we will call them $\Omega\Psi$ rules,

\begin{align}
&\label{Omega1}(\forall x, y\in\Omega)( x+y\in\Omega, xy\in\Omega),\\
&\label{Omega2}(\forall x\in\Omega)( \forall y\in\mathbb A)( xy\in\Omega),\\
&\label{Omega3}(\forall x)(x\in \mathbb A\Rightarrow  x^{-1}\in\mathbb A),\\
&\label{Omega4}(\forall x\neq 0)(x\in\Omega  \Leftrightarrow \ x^{-1}\in\Psi),\\
&\label{Omega5}(\forall x\in \Psi)(\forall y\in \A)(\exists z\in\Omega) (xz=y).
\end{align}

Referring to the set $\Omega$, an equivalence relation is defined by
\[x\approx y \Leftrightarrow x-y \in\Omega,\ \ \ \mbox{for}\ x,y\in\F.\]

We say that $x$ \emph{is infinitely close to} $y$, when the relation $x\approx y$ holds.

  We will show in section 6 below that in (Euler 1755), Euler explicitly discusses rules equivalent to $\Omega\Psi$ rules.
  
\subsection{Archimedean axiom} Below we present some equivalent forms of the Archimedean axiom.
\begin{enumerate}\itemsep 0mm
  \item [(A1)] $(\forall x,y\in \F)(\exists n\in\N)(0<x<y\Rightarrow nx>y)$,
  \item [(A2)] $(\forall x\in\mathbb F)(\exists n\in\N)(n>x) $,
\item[(A3)] $\Omega=\{0\}$.
\end{enumerate}

 A1 and A2 are well-known versions both in mathematical and historical contexts. Due to Euler's development, we prefer the A3 version.  It reads that in a non-Archimedean field, the set of infinitesimals contains at least one non-zero element, say, $\varepsilon$.      Then, $\tfrac \varepsilon n$, as well as $n\varepsilon$   are infinitesimals too. Furthermore,  by $\Omega\Psi$ rules,   reciprocals for these elements, i.e.  $\tfrac n \varepsilon$, $\tfrac 1 {n\varepsilon}$, are  infinitely large numbers.

In what follows, let $\Omega_0$ stand for non-zero infinitesimals, i.e.
  \[\Omega_0=\Omega\setminus\{0\}.\]

\subsection{Real numbers}
 The field of real numbers   is  an  ordered field in which
 every Dedekind cut $(L,U)$ of $(\F,<)$ satisfies the following condition,  the so-called completeness  axiom, CA in short,
\[(\exists x\in\F)(\forall y\in L)(\forall z\in U)(y\leq x\leq z).\]

The above definition applies the theorem that every two ordered fields satisfying CA are isomorphic. In this sense, the field of real numbers is the unique complete ordered field. Moreover, any Archimedean field is isomorphic to a subfield of real numbers. As a result, any field extension of real numbers is non-Archimedean and includes infinitely small and infinite numbers.

\subsection{Hypperreals}
 
We define the set of  hypperreals  as the quotient set on the set of all sequences of real numbers, $\R^{\N}$, with respect to a specific relation defined on the set of indexes $\N$. To this end, 
  we need a notion of ultrafilter on $\N$.

 A family of sets $\mathcal U\subset \mathcal P (\N)$ is an ultrafilter on $\N$ iff (1) $\emptyset \notin \mathcal U$, (2) if $A,B\in\mathcal U$, then $A\cap B\in\mathcal U$, (3) if $A\in\mathcal U$ and $A\subset B$, then $B\in\mathcal U$, (4) for each  $A\subset \N$, either $A$ or its complement $\N\setminus A$ belongs to $\mathcal U$.

  The family of  sets with  finite complements satisfies conditions (1)--(3) listed in the definition of an ultrafilter. By  Zorn's lemma, this family  can be  extended to an  ultrafilter.   Let $\mathcal U$ denote a fixed ultrafilter  on $\N$ containing  every  subset with a finite complement.

Since $\N$ and sets of the form $\{k, k+1, k+2,...\}$, in short $\{n\in\N: n\geq k\}$, for every $k\in\N$, comply 
with the proviso  ``have finite complements", they all belong to the ultrafilter $\mathcal U$,
\[\N\in \mathcal U, \ \ \{n\in\N: n\geq k\}\in \mathcal U, \ \ \mbox{for\ every}\ k.\] 

These sets will do to check examples concerning infinitesimals and infinite numbers, which we present below.

 In the set $\R^{\N}$ we define an equivalence relation by
\[(r_n)\equiv (s_n)\Leftrightarrow \{n\in \N:\ r_n=s_n\}\in\mathcal U.\]

Let $\R^*$  denote  the reduced product $\R^{\N}/_\equiv$.

Clearly, the equality of hyperreals is defined by 
\[[(r_n)]= [(s_n)]\Leftrightarrow \{n\in \N:\ r_n=s_n\}\in\mathcal U.\]

Since $\mathcal U$ is the ultrafiter, we also get  the following condition
\[[(r_n)]\neq [(s_n)]\Leftrightarrow \{n\in \N:\ r_n\neq s_n\}\in\mathcal U.\]

It is used to determine the multiplicative inverse of a hyperreal, namely
\[[(r_n)]^{-1}=[(r_1^{-1},r_2^{-1}, r_3^{-1},...)].\]

To illustrate how it works, an inequality 
\[[(0,0,0,...)]\neq [(r_1,r_2, r_3,...)]\]

translates into condition $\{n\in\N: r_n\neq 0\}\in\mathcal U$, rather than $\{n\in\N: r_n\neq 0\}=\N$.

 New sums and products are defined pointwise, that is
\[[(r_n)]+[(s_n)]=[(r_n+s_n)],\ \ \  [(r_n)]\cdot [(s_n)]=[(r_n\cdot s_n)].\]

A total order on $\R^*$ is  defined by
\[[(r_n)]<[(s_n)]\Leftrightarrow \{n\in \N :r_n<s_n\}\in\mathcal U.\footnote{In these definitions, we adopt a standard convention to use the same signs for the relations and operations  on $\R$ and $\R^*$.}\]

A standard  real number $r\in \R$  is represented by the class $[(r,r,r,...)]$. In what follows, we will use 
 the letter  $r$ for    hyperreal  number $[(r,r,r,...)]$.

By a straightforward checking, we find $(\R^*,+,\cdot,0,1,<)$ is an ordered field.

We also introduce the set of non-standard natural numbers $\N^*$, namely
\[\N^*=\{[(n_j)]\in\R^*: (n_j)\subset \N^{\N}\}.\] 

\subsection{Infinitesimals, infinite numbers, finite numbers}

The field of hyperreals is non-Archimedean and enables one, in a sense,  to \textit{touch} infinitesimals and infinite numbers. 
To demonstrate that hyperreals extend the field of reals, we require a single non-zero infinitesimal.
The hyperreal number $\varepsilon=[(1,\tfrac 12, \tfrac 13,...)]$ is  positive and smaller than every  positive real number, i.e.
\[[(1,\tfrac 12, \tfrac 13,...)]<[(r,r,r,...)],\ \ \mbox{for\ every}\ \ r\in\R_+.\]

For a proof, take $r>0$. In real analysis, $\lim\limits_{n\rightarrow \infty}1/n=0$.\footnote{Indeed, it is yet another form of the Archimedean axiom.} 
It means, there is  an index $k$ such that $1/n<r$, for all $n>k$. In terms of the elements of the ultrafiter, $\{n\in\N: n> k\}\in \mathcal U$, or 
\[\{n\in\N: \tfrac 1n<r\}\in\mathcal U.\]

Due to the definition of the total order in $\R^*$, it means
\[[(1,\tfrac 12, \tfrac 13,...)]<[(r,r,r,...)].\]

Thus $\varepsilon\in\Omega$. The hyperreal 
$$\varepsilon^2=[(1,\tfrac 1{2^2}, \tfrac 1{3^2},...)]$$ 
is another infinitesimal. Generally, if $(r_n)$ is a null-sequence, then $[(r_n)]$ is infinitesimal.

Since $\varepsilon\in\Omega$, then $\varepsilon^{-1}=[(1,2,3,...)]$ is an infinitely large number. Indeed,  $N=[(1,2,3,..)]$ exemplifies   infinite number. Other infinite numbers are 
$$N+1=[(2,3,4,...)],\ \ N^2=[(1,4,9,...)],\ \  N!=[(1!, 2!, 3!,..)].$$

\subsection{Binomial coefficients}

Here is how we define the binomial coefficients $\binom Nk$ for infinite $N$, $k$, and $N,k\in\N^*$. Let us start with the case
$N=[(n_1,n_2,n_3,...)]$ is infinite, $k$ is finite. We put
\[\binom Nk=\Big[\Big(\binom {n_1}k, \binom {n_2}k, \binom {n_3}k,...\Big)\Big].         \]

When $K=[(k_1,k_2,k_3,...)]$ is an infinite number, we put
\[\binom NK=\Big[\Big(\binom {n_1}{k_1}, \binom {n_2}{k_2}, \binom {n_3}{k_3},...\Big)\Big].         \]

These new binomial coefficients meet the standard identities, like $\binom nk= \binom n{n-k}$.

For infinite $N$ and finite or infinite $k$ the implicit assumptions of Euler's arguments obtains in the following form
\[ \dfrac {N-k}{kN}\approx \dfrac{1}{k},\ \ \ \ \ \binom Nk \dfrac{1}{N^k} \approx \dfrac 1 {k!}.\]
\subsection{*Maps} 
Let $f$ be a real map, i.e. $f:\R\mapsto \R$. By  $f^*$ we mean a map
  $f^*:\R^*\mapsto \R^*$ defined by
\[f^*([(r_n)]) =[(f(r_1),f(r_2),...)].\]

If $r\in \R$, then $f^*(r)=[(f(r),f(r),...)]$.
  Since we identify real number $r$ with hyperreal $[(r,r,...)]$,  the equality $f^*([(r,r,...)])=f(r)$ obtains.

Under this definition, 
\[\sin^*[(r_n)]=[(\sin r_1,\sin r_2,...)], \ \ \cos^*[(r_n)]=[(\cos r_1,\cos r_2,...)].\]

Since for every $n$ the identity $\sin^2 r_n+\cos^2 r_n=1$  holds, we have 
\[(\sin^*x)^2+(\cos^* x)^2=1.\]

Similarly, every trigonometric identity can be transferred into an identity involving the maps $\sin^*$ and  $\cos^*$.

The above definition applies also to the exponential map $a^x$. Thus, when $a\in \R_+$, then
\[a^{[(r_n)]}=[(a^{r_1}, a^{r_2},...)].\]

When $a=[(a_1, a_2,...)]$ is hyperreal, then  
\[a^{[(r_n)]}=[({a_1}^{r_1}, {a_2}^{r_2},...)].\]

Then, for $a, b, x, y \in \R^*_+$, we have $(ab)^x=a^xb^x$, $a^xa^y=a^{x+y}$, etc.

Let us go back to Euler's supposition  $\sin z\approx z$, given $z$ is infinitesimal. Clearly, $|\sin^* z|\leq 1$, i.e.,  $\sin^* z$ is assignable or infinitesimal. To get a contradiction, suppose $\sin^* z$ is an assignable number for some infinitesimal $z$.  Since the real sine  takes  the segment $[0,\tfrac \pi 2]$ onto the segment of real 
numbers $[0,1]$, the $\sin^*$ is also a surjective map.  As a result, for some assignable number 
$x$, the equality obtains $\sin^* z=\sin^*x$. On the segment $[0,\tfrac \pi 2]$, the real sine is one-to-one,   therefore 
$\sin^*$ is one-to-one on the segment of hyperreals $[0,\tfrac \pi 2]$. Hence, $z=x$, an infinitesimal is equal to an assignable number. This is impossible. 

Another way to reach that result is to refer to the continuity of the real map sine. Its continuity at the point $0$ translates into the condition $\sin z\approx 0$, given $z\approx 0$. Generally, one has to assume some characteristics of the real map $f$ to make some claims on $f^*$. Euler, however, did not define a sine. Thus, we are trapped and have to speculate on what grounds $\sin z\approx z$. 

The assumption $\cos^* z\approx 1$, given $z\approx 0$, follows form the Pythagorean identity for 
$\sin^*$ and $\cos^*$, and the just discussed result $\sin^*z\approx 0$.

\subsection{Hyperfinite sums}

Let $N=[(n_1,n_2,n_3,...)]$ be an infinite number, $N\in\N^*$, let $a$ be a real number. We set
\[\sum_{j=0}^{N} a^j=\Big[\Big(\sum_{i=0}^{n_1} {a}^i,\sum_{i=0}^{n_2} {a}^i,\sum_{i=0}^{n_3} {a}^i,... \Big)\Big].\]

Thus
\[1+a^1+a^2+a^3+...+a^N\]
is another symbol for the hyperreal number
\[\Big[\Big(\sum_{i=0}^{n_1} {a}^i,\sum_{i=0}^{n_2} {a}^i,\sum_{i=0}^{n_3} {a}^i,... \Big)\Big].\]


Since for every finite $n$ the equality obtains
$$\sum_{i=0}^{n} {a}^i = \frac{1- a^{n+1}}{1-a},$$

we easily get 
\begin{align*}
\sum_{j=0}^{N} a^j &=\Big[\Big(\frac{1- a^{n_1+1}}{1-a}, \frac{1- a^{n_2+1}}{1-a}, \frac{1- a^{n_3+1}}{1-a}, ...\Big)\Big]\\
&= \frac{1- a^{N+1}}{1-a}.
\end{align*}

When $a$ is a hyperreal number, $a=[(a_1,a_2,a_3,...)]$, by putting 
\[\sum_{j=0}^{N} a^j =\Big[\Big(\sum_{i=0}^{n_1} {a_1^i},\sum_{i=0}^{n_2} {a_2^i},\sum_{i=0}^{n_3} {a_3^i},... \Big)\Big],\]
we reach the same result, namely 
\[\sum_{j=0}^{N} a^j=\frac{1- a^{N+1}}{1-a}.\]

\section{Infinitesimals and infinite numbers in \textit{Institutiones calculi differentialis}}

Chapter III (\S\S 74--111) of \textit{Institutiones calculi differentialis} (Euler 1755) begins with philosophical considerations like whether matter consists of indivisible or infinitely divisible parts, or whether supposed ultimate parts of matter are extended or not. Euler finds these debates inconclusive and turns to mathematics by saying:  

\begin{quote}Even if someone denies that infinite numbers really exist in this world, still in mathematical speculations there arise questions to which answers cannot be given unless we admit an infinite number. (Euler 2000, 50)\end{quote}

\subsection{Three kinds of quantities. Infinite numbers}

In section 82, Euler offers a definition of infinite numbers:

\begin{quote}[...] this quantity is so large that it is greater than any
finite quantity and cannot not be infinite. To designate a quantity of this
kind we use the symbol $\infty$, by which we mean {a quantity greater than any
finite or assignable quantity}. (Euler 2000, 50)\end{quote}

\pv  In symbols, given $\A$ stands for assignable numbers:
\begin{equation}N\ \mbox{is\ infinite\ number}\Leftrightarrow (\forall a\in\mathbb A)(N>a).\end{equation}

The word \emph{quantitas}  (quantity) is reiterated over and over again throughout the discussed chapter. It refers to infinite, assignable, and infinitesimals numbers. \textit{Quantitas}  is Latin translation of the Greek term  \foreignlanguage{polutonikogreek}{m'egejos} (magnitude). In Euclid's \textit{Elements},  this general term   covers line segments, triangles, convex polygons, circles,  angles, arcs of circles, and solids.  Magnitudes of the same kind
(line segments being of one kind, triangles of another, etc.) are compared in terms of \textit{greater}-- \textit{lesser}}.  Transitivity and the law of trichotomy constitute the mathematical sense of \textit{greater than} relation.  Greeks also considered addition and subtraction (lesser from greater) of magnitudes. The additive and ordinal structure of magnitudes enabled Euclid to introduce the proportion theory.   Book V of the \textit{Elements} develops it systematically. 

 It was the basic ancient Greek technique regarding comparing magnitudes. The Archimedean axiom was an essential part of that theory and significantly restricted the concept of magnitude.  Descartes (Descartes 1637) introduced a  new operation on line segments: the product. He employed it in such a way that resulting arithmetic satisfied the rules of an ordered field.  While processing formulas, he substitutes equality of quotients for proportions: instead of proportions of line segments, such as $a:b::c:d$, he puts $\tfrac ab=\tfrac cd$.\footnote{See (B{\l}aszczyk, Petiurenko 2019), (B{\l}aszczyk 2022).} As a result, the Archimedean axiom lost its gravity and for a long time has not been discussed. Only (Stolz 1885) re-introduced it to modern mathematics and re-established its role in foundational studies.

Euler, as we have already seen, processes formulas according to the rules of an ordered field. 
The above genealogy provides a rationale for his manner of naming quotients ``geometrical ratios" (see the next sections).  He even applies the symbol $a:b$, used for the proportion of magnitudes through the 17th and 18th centuries, though he processes these objects as actual fractions $\tfrac ab$.

The term Archimedean axiom never occurs in Euler's writings. Yet, due to the explicit negation of the axiom A3,  his implicit ordered field is non-Archimedean.

\subsection {Infinitesimals}

Section 83  provides a definition of  an infinitesimal number:

\begin{quote}There is also a definition of {the
infinitely small quantity as that which is less than any assignable quantity}.
If a quantity is so small that it is less than any assignable quantity, then
it cannot not be $0$, since unless it is equal to $0$ a quantity can be assigned
equal to it, and this contradicts our hypothesis. To anyone who asks what
an infinitely small quantity in mathematics is, we can respond that it really
is equal to 0. (Euler 2000, 51)\end{quote}

\pv In symbols
\begin{equation}\varepsilon\ \mbox{is\ infinitesimal}\Leftrightarrow (\forall a\in\mathbb A)(0<|\varepsilon|<|a|).\end{equation}

Clearly, the above formalization does not comply with Euler's words ``cannot not be $0$". However, his concept of equality is ambiguous: in some contexts, it is the standard, say, strict equality, in others it also means \textit{is infinitely close}. The latter term comes from the non-standard analysis. We substantiate it in sections that follow.

\subsection{Two ways of comparing  zeros}

Euler explicitly claims that in the analysis, the other name for zero is infinitely small.

\begin{quote} If we accept the notation used in the analysis of the infinite, then $dx$ indicates the quantity that is infinitely small,
so that both $dx=0$ and
$a dx = 0$, where $a$ is any finite quantity.  Despite this, the geometric ratio
$adx : dx$ is finite, namely $a : 1$. For this reason these two infinitely small
quantities $dx$ and $adx$, both being equal to $0$, cannot be confused when
we consider their ratio. (Euler 2000, 51--52)\end{quote}

\pv If we adopt the interpretation of the formula   $dx=0$ as $dx\approx 0$, Euler's two ways of comparing zeros get an obvious meaning. Infinitesimals, say $\epsilon,\delta$, can be compared as greater-lesser by $\epsilon - \delta$  or $\dfrac \epsilon\delta$. The first term is always infinitesimal, while the second an infinitesimal, finite, or infinite. 
The  following  passage clarifies his view:

\begin{quote}[1] Although {two zeros are equal to each
other}, so that there is no difference between them, nevertheless, since we
have {two ways to compare them}, either arithmetic or geometric, let us look
at quotients of quantities to be compared in order to see the difference.
[2] The {arithmetic ratio} between any two zeros is an equality.   This is not the
case with a {geometric ratio}.   [3] We can easily see this from this geometric proportion
$2:1=0 : 0$, in which the fourth term is equal to $0$, as is the third.
From the nature of the proportion, since the first term is twice the second, it is necessary
that the third is twice the fourth. (Euler 2000, 51; numerals in square brackets added).\end{quote}

\pv To elaborate. Ad 1.  The   ``equality of zeros"  means the infinitely close relation, $\epsilon\approx\delta$.  Ad 2. ``Arithmetic ratio" means  $\epsilon-\delta$. Ad 3. The ``geometric ratio" means  the quotient $\frac \epsilon\delta$. Substituting $\epsilon$ for $0$ in Euler's formula $2:1=0:0$, we  get $2:1=\frac {2\epsilon}\epsilon$.

In Greek mathematics, ratio, like $a:b$, made sense only in a proportion, say, \mbox{$a:b::c:d$.} Euler, same as Newton, defines number as a ratio. Consequently, terms such as $adx:dx$ are numbers.  The phrase ``geometric ratio"  refers back to Euclid's proportion of magnitudes. However, considering actual practice, Euler processes numbers according to the rules of an ordered field, rather than propositions of Euclid's Book V.

\subsection{Arithmetic of infinitesimals}

\begin{quote}Since the infinitely small is actually nothing, it is clear that a finite
quantity can neither be increased nor decreased by adding or subtracting an
infinitely small quantity. Let $a$ be a finite quantity and let $dx$ be infinitely
small.  Then $a+dx$ and $a-dx$, or, more generally, $a\pm n dx$, are equal to $a$. 
[...] On the other hand, the  geometric ratio is clearly of equals, since $\frac{a\pm n dx}{a}=  1$. (Euler 2000, 52)\end{quote}

Given $a$ is a finite number, Euler's claims obviously translate into
\[a\pm ndx\approx a, \ \ \ \frac{a\pm n dx}{a}\approx  1.\]

Regarding products of infinitesimals, Euler writes

\begin{quote} Since the infinitely small quantity $dx$ is actually equal to 0, its square
$dx^2$, cube $dx^3$, and any other $dx^n$, where n is a positive exponent, will
be equal to 0, and hence in comparison to a finite quantity will vanish.
However, even the infinitely small quantity $dx^2$ will vanish when compared
to $dx$. The ratio of $dx\pm dx^2$ to $dx$ is that of equals, whether the comparison
is arithmetic or geometric. There is no doubt about the arithmetic; in the
geometric comparison $dx\pm dx^2 : dx = \dfrac{dx\pm dx^2}{dx}= 1\pm dx = 1$. (Euler 2000, 52)\end{quote}

In our interpretation
\[dx^n\approx dx,\ \ \  \frac{dx+dx^n}{dx}\approx 1.\]

Euler also considers infinitesimals of the form $\sqrt{dx}$ and claims that $\sqrt{dx}=dx$. In our interpretation, it means $\sqrt{dx}\approx dx$.

 \subsection{Infinitesimals and infinite numbers} 
\begin{quote}[1] It should
be noted that the fraction $1/z$ becomes greater the smaller the denominator
$z$ becomes. [2] Hence, if $z$ becomes a quantity less than any assignable quantity,
that is, infinitely small, then it is necessary that the value of the fraction $1/z$
becomes greater than any assignable quantity and hence infinite. For this
reason, if $1$ or any other finite quantity is divided by something infinitely
small or $0$, the quotient will be infinitely large, and thus an infinite quantity.
Since the symbol  $\infty$  stands for an infinitely large quantity,  we have the
equation $\frac a{dx}=\infty$. [3] The truth of this is clear also when we invert
$\frac{a}{\infty}=dx=0$. (Euler 2000, 53; numerals in square brackets added)\end{quote}

\pv The first sentence of this passage explains the currently obvious rule of ordered fields, namely:
 \[0<z<x\Rightarrow \frac 1x<\frac 1z.\]

The second justifies  our (\ref{Omega4}) $\Omega\Psi$ rule, 
\[z\in\Omega_0\Rightarrow \frac 1z\in\Psi.\]

The third derives the reverse implication,
\[N\in\Psi,\ a\in\mathbb A \Rightarrow \frac a N\in\Omega_0.\]

Then Euler explains why there are infinitely many infinite numbers: 

\begin{quote}Since $a/dx$ is an infinite quantity $A$, it is clear that the quantity $A/dx$
will be a quantity infinitely greater than the quantity $A$. This can be seen
from the proportion $a/dx : A/dx = a : A$, that is, as a finite number to
one infinitely large. There are relations of this kind between infinitely large
quantities, so that some can be infinitely greater than others. Thus, $a/dx^2$
is a quantity infinitely greater than $a/dx$; if we let $a/dx = A$, then $a/dx^2 =
A/dx$. In a similar way $a/dx^3$ is an infinite quantity infinitely greater than
$a/dx^2$, and so is infinitely greater than $a/dx$. {We have, therefore an infinity
of grades of infinity, of which each is infinitely greater than its predecessor}.
If the number m is just a little bit greater than n, then $a/dx^m$ is an infinite
quantity infinitely greater than the infinite quantity $a/dx^n$. (Euler 2000, 55)\end{quote}

\pv In our interpretation, both $1/\epsilon$ and $1/\epsilon^2$ belong to $\Psi$, given $\epsilon\in\Omega_0$.
Yet, comparing them in terms of ``geometric ratio", that is, $\frac{1/\epsilon^2}{1/\epsilon}$,  we observe the second is greater than the first one.

\subsection{Assignable numbers and  $\Omega\Psi$ products} 
In nonstandard analysis, finite numbers include assignable numbers, $\A\subset \mathbb L$. Euler does not consider infinitesimals to be finite. That is why, in our interpretation, his assignable numbers form the set $\mathbb L\setminus \Omega$. Indeed, $1, 2$, and other integers exemplify assignable numbers. Since $a+\epsilon\approx a$, both standard number $1$ and $1+\epsilon$ are assignable.

\begin{quote}Now that we have been warned about the grades of infinities, we will
soon see that it is possible not only for the product of an infinitely large
quantity and an infinitely small quantity to produce a finite quantity, as
we have already seen, but also that a product of this kind can also be
either infinitely large or infinitely small. (Euler 2000, 55--56)\end{quote}

The examples that follow, we present  this in non-standard stylization. Given $a\in\A$, $\epsilon\in\Omega$, so $a/\epsilon$ and  $a/\epsilon^2$ are infinite. Then
\[\frac a\epsilon \epsilon\in\A,\ \ \ \frac a\epsilon \epsilon^2\in\Omega,\ \ \ \frac{a}{\epsilon^2}\epsilon\in\Psi. \] 

The remaining sections of the discussed chapter consider infinite sums.   Therein Euler introduces objects such as $1+x+x^2+...+x^\infty$, which we studied above in terms of hyperfinite sums. 

\section{Summary}
Euler's technique of determining $\sin v$, for $v=\frac mn\frac \pi 2$ consists in summing up first terms of the series 
\[v-\frac{v^3}{3!}+\frac{v^5}{5!}-\frac{v^7}{7!}+...\,,\] 
taking into account numerical approximation of $\pi$. To derive that series, he takes an infinite number $j$ and  an infinitesimal 
$\omega$ such that $v=j\omega$. By transforming the  identity
\[\sin j\omega=\dfrac{(\cos\omega+i\sin\omega)^j+ (\cos\omega-i\sin\omega)^j}{2i}\]
according to the binomial theorem, and replacing  $\sin \omega$ with $\omega$,  $\cos \omega$ with $1$, and $\binom jk \tfrac 1{j^k}$ with $\tfrac 1{k!}$,  he gets the currently well-know series for sine. Series for cosine Euler derives through the same technique.
 
 We showed that most of Euler's argument is valid in any non-Archimedean field. Analyzing his  \textit{Institutiones calculi differentialis}, we proved he explicitly examined the rules of such a  field.  Finally, we presented hyperreal numbers, which enable one to justify all the assumptions of Euler's proof.

\begin{acknowledgement}

We thank an anonymous reviewer for insightful suggestions.

\pv
The first author is supported by the National Science Centre, Poland grant 2018/31/B/HS1/03896  Infinite and Infinitesimals.

The second author is supported by the National Science Centre, Poland grant 2018/31/B/HS1/03896  Infinite and Infinitesimals.

\end{acknowledgement}

%
%
%

\end{document}